\documentclass[12pt]{article}
\usepackage{amsfonts}
\usepackage[vcentermath]{youngtab}

\title{Betti numbers of Springer fibers in type $A$}
\author{Lucas Fresse}

\newcounter{numpartie}
\newcounter{numparagraphe}[numpartie]

\renewcommand{\thenumparagraphe}{\arabic{numpartie}.\arabic{numparagraphe}}

\newcommand{\partie}[1]{\refstepcounter{numpartie}
\subsection*{\thenumpartie. #1}}
\newcommand{\paragraphe}{\refstepcounter{numparagraphe}
\paragraph{\,}\ \textit{\thenumparagraphe. }}

\newcommand{\C}{\mathbb{C}}
\newcommand{\N}{\mathbb{N}}
\newcommand{\Q}{\mathbb{Q}}

\newcommand{\F}{F}
\newcommand{\B}{{\mathcal F}}

\newcommand{\nilp}{u}

\newcommand{\Sym}{S}

\newcommand{\cqfd}{\hfill$\sqcup\!\!\!\!\sqcap$}

\begin{document}

\maketitle

\begin{center}
{\em
Department of Mathematics,
the Weizmann Institute of Science,
Rehovot 76100,
Israel}\\
E-mail: lucas.fresse@weizmann.ac.il
\end{center}

\begin{abstract}
We determine the Betti numbers of the Springer fibers in type $A$.
To do this, we construct a cell decomposition of the Springer fibers.
The codimension of the cells is given by an analogue of the Coxeter length.
This makes our cell decomposition well suited for the calculation of Betti numbers.
\end{abstract}


\partie{Introduction}

\paragraphe 
Let $V$ be a $\C$-vector space of dimension $n\geq 0$
and let $\nilp:V\rightarrow V$ be a nilpotent endomorphism.
We denote by $\B$ the (algebraic) variety of complete flags of $V$
and by $\B_\nilp$
the subset of $\nilp$-stable complete flags, i.e.
flags $(V_0,...,V_n)$ such that $\nilp(V_i)\subset V_i$
for all $i$.
The variety $\B$ is projective,
and $\B_\nilp$ is a projective subvariety of it.
The variety $\B_\nilp$ is called {\em Springer fiber}
since it can be seen as the fiber over $\nilp$ of the Springer resolution
of singularities of the cone of nilpotent endomorphisms of $V$
(see for example \cite{Slodowy}).

Springer constructed representations of the
symmetric group $\Sym_n$ on the cohomology spaces $H^{*}(\B_\nilp,\Q)$  (see \cite{Springer}). The characters of these
representations were determined by Lusztig in \cite{Lusztig}.
More explicitly he connected the multiplicities
of irreducible summands of $H^{*}(\B_\nilp,\Q)$
with the coefficients of Kostka polynomials.
This allows
to calculate the Betti numbers $b_m=\mathrm{dim}\,H^{2m}(\B_\nilp,\Q)$. The aim of this article is to give a more direct
calculation of them.

\paragraphe
Following \cite{DeConcini}, a finite partition of a variety $X$ 
is said to be an $\alpha$-partition if the subsets in the partition
can be indexed $X_1,...,X_k$ so that $X_1\cup ...\cup X_l$
is closed in $X$ for $l=1,2,...,k$.
Thus each subset in the partition is a locally closed subvariety of $X$.
An $\alpha$-partition into subsets which are isomorphic to affine spaces is called 
a cell decomposition.
If $X$ is a projective variety with a cell decomposition, then 
the cohomology of $X$ vanishes in odd degrees and
$\mathrm{dim}\,H^{2m}(X,\Q)$ is the number of $m$-dimensional cells
(see \ref{cohomology}).

It is known from \cite{Shimomura} and \cite{Spaltenstein}
that $\B_\nilp$ admits a cell decomposition,
and there are also many references 
%giving a construction 
proving the existence
for other types
(see \cite{Spaltenstein2} or \cite{NanhuaXi}) or
more general contexts (Springer fibers of any type in \cite{DeConcini},
partial $\nilp$-stable flags in \cite{Nakajima}).
A simple manner to construct a cell decomposition of $\B_\nilp$
is to take the intersection with the Schubert cells of the flag
variety, then we obtain a cell decomposition
provided that the Schubert cells are defined according to appropriate conventions
(see \cite{Shimomura} or \ref{section-intersection-Schubert}).
However the dimension of the cells
is given by a complicated formula, 
it makes this cell decomposition not practical to compute Betti numbers.
We construct a different cell decomposition which is better suited for the
calculation of Betti numbers.

\paragraphe
\label{section-theoreme-principal}
Let $\lambda(\nilp)=(\lambda_1\geq ...\geq \lambda_r)$ be the lengths of the Jordan
blocks of  $\nilp$ and let $Y(\nilp)$ be the Young diagram of rows
of these lengths.
If $\mu_1,...,\mu_s$ are the lengths of the columns of $Y(\nilp)$,
recall from \cite[\S II.5.5]{Spaltenstein} that
$\mathrm{dim}\,\B_\nilp=\sum_{q=1}^s{\mu_q(\mu_q-1)}/2$.

A standard tableau of shape $Y(\nilp)$ is a numbering of the boxes of $Y(\nilp)$ by $1,...,n$
such that numbers in the rows increase to the right and
numbers in the columns increase to the bottom.

We call {\em row-standard tableau of shape $Y(\nilp)$}
%$\lambda$}
a numbering of the boxes of $Y(\nilp)$ by $1,...,n$ 
such that numbers in the rows increase to the right.
%We call nonordered partition a sequence of nonnegative integers of finite sum.
Let $\tau$ be a row-standard tableau.
We call {\em inversion} %of $\tau$ 
a pair of numbers $i<j$
in the same column of $\tau$ and such that one of the following conditions
is satisfied: \\
\ \ - \ $i$ or $j$ has no box on its right
and $i$ is below $j$, \\
\ \ - \ $i,j$ have respective entries $i',j'$ on their right, and $i'>j'$.

\smallskip
\noindent
For example $\tau=\mbox{\scriptsize $\young(248,367,15)$}$ has four inversions: 
the pairs $(1,2)$, $(4,6)$, $(5,6)$, $(7,8)$. \\
Let $n_\mathrm{inv}(\tau)$ be the number of inversions of $\tau$.
We see that $n_\mathrm{inv}(\tau)=0$ if and only if $\tau$ is a standard tableau.
For $\nilp=0$ the diagram $Y(\nilp)$ has only one column, hence
$\tau$ is equivalent to a permutation
($\sigma\in S_n$ corresponds to the tableau numbered by $\sigma_1,...,\sigma_n$
from top to bottom) and $n_\mathrm{inv}(\tau)$ is the usual
inversion number for permutations. 

Our main result is the following

\medskip
\noindent
{\bf Theorem}
{\em
%(a) 
The variety $\B_\nilp$ has a cell decomposition $\B_\nilp=\bigcup_\tau C(\tau)$
parameterized by the row-standard tableaux of shape $Y(\nilp)$,
and such that $\mathrm{codim}_{\B_\nilp}\,C(\tau)=n_\mathrm{inv}(\tau)$. %\\
%(b)
%Let $d=\mathrm{dim}\,\B_N$. 
%For $m\geq 0$, the Betti number $b_m$
%is the number of row-standard tableaux 
%$\tau$
%such that $n_\mathrm{inv}(\tau)=d-m$.
}

\medskip

And we deduce:

\medskip
\noindent
{\bf Corollary}
{\em
Let $d=\mathrm{dim}\,\B_\nilp$. 
For $m\geq 0$, the Betti number $b_m:=\mathrm{dim}\,H^{2m}(\B_\nilp,\Q)$
is the number of row-standard tableaux 
$\tau$ of shape $Y(\nilp)$
such that $n_\mathrm{inv}(\tau)=d-m$.
}

\medskip
If $\nilp=0$, then $\B_\nilp$ is the whole flag variety $\B$,
and we get the classical formula giving the Betti numbers of the flag variety.
In general, we find that the dimension of the cohomology space
of maximal degree is the number of standard tableaux of shape $Y(\nilp)$.
This is also classical, since
the Springer representation of $S_n$ on the cohomology in
maximal degree is irreducible and isomorphic to the Specht module corresponding
to the Young diagram $Y(\nilp)$, whose dimension is precisely 
the number of standard tableaux of shape $Y(\nilp)$
(see \cite{Springer}).
We also recall in \ref{paragraphe-BuT} that $\B_\nilp$ is equidimensional and that its
components are parameterized by standard tableaux.

\paragraphe
\label{standardization}
Let us make more precise the relation
between standard and row-standard tableaux.
If $T$ is standard, then the shape of its subtableau $T[1,...,i]$ of entries $1,...,i$
is a subdiagram $Y_i(T)\subset Y(\nilp)$.
In that way a standard tableau $T$ is equivalent to the data of a complete chain of subdiagrams
$\emptyset=Y_0(T)\subset Y_1(T)\subset ...\subset Y_n(T)=Y(\nilp)$.
We call (ordered) partition of $n$ a decreasing sequence of nonnegative integers whose sum is $n$.
The lengths of the rows of $Y_i(T)$ form a partition of $i$ \
$\lambda^{(i)}=(\lambda^{(i)}_1\geq ...\geq \lambda^{(i)}_r)$. % of sum $i$
%and, similarly, $T$ defines a maximal chain of partitions
%$\emptyset=\lambda^{(0)}\subset \lambda^{(1)}\subset ...\subset \lambda^{(n)}=\lambda(\nilp)$
%(we denote by $\subset$ the partial inclusion order on partitions, which means $\lambda^{(i)}_p\leq
%\lambda^{(i+1)}_p$ for every $p$).
%In that way a standard tableau is equivalent to the data of a maximal chain from $\emptyset=(0,0,...,)$
%to $\lambda(\nilp)$ in the poset of partitions.
In that way, $T$ is also equivalent to a maximal chain of partitions
$\emptyset=\lambda^{(0)}\subset \lambda^{(1)}\subset ...\subset \lambda^{(n)}=\lambda(\nilp)$
from $\emptyset=(0,0,...)$ to $\lambda(\nilp)$
(we denote by $\subset$ the partial inclusion order on partitions, which means $\lambda^{(i)}_p\leq
\lambda^{(i+1)}_p$ for every $p$).

If $\tau$ is row-standard, then
%the shape of its subtableaux are not subdiagrams (the length of rows does not necessarly decrease),
the lengths of the rows of its subtableau of entries $1,...,i$ form a 
sequence of nonnegative integers
$\pi^{(i)}=(\pi^{(i)}_1, ..., \pi^{(i)}_r)$ of sum $i$
(not necessarily weakly decreasing).
In that way $\tau$ is equivalent to the data of a maximal chain
of finite sequences of nonnegative integers
$\emptyset=\pi^{(0)}\subset \pi^{(1)}\subset ...\subset \pi^{(n)}=\lambda(\nilp)$
from $\emptyset$ to $\lambda(\nilp)$%, in the poset
\ (where %the partial order 
$\subset$ means
$\pi^{(i)}_p\leq\pi^{(i+1)}_p$ for every $p$).

If we arrange the entries in each column of $\tau$ in increasing order
to the bottom, then we get a standard tableau
that we denote by $\mathrm{st}(\tau)$.
We will call it the standardization of $\tau$.
Similarly if we arrange the terms of each sequence $\pi^{(i)}$
in decreasing order, then we get a partition $\mathrm{ord}(\pi^{(i)})$,
the partitions $(\mathrm{ord}(\pi^{(i)}))_{i=0,...,n}$ form a maximal chain
from $\emptyset$ to $\lambda(\nilp)$
and $\mathrm{st}(T)$ is the standard tableau which corresponds
to it.
As we show in \ref{proposition-coxeter}, the inversion number of $\tau$ can be interpreted
as a minimal number of elementary operations which allow to
transform $\tau$ into its standardization $\mathrm{st}(\tau)$.

%\smallskip

%We also need the dominance order on diagrams and standard tableaux.
%If $Y$ and $Y'$ are two Young diagrams
%of rows of (possibly zero) lengths $(\lambda_1\geq...\geq \lambda_r)$
%and $(\lambda'_1\geq...\geq \lambda'_r)$,
%we write $Y\preceq Y'$ if $\lambda_1+...+\lambda_p\leq \lambda'_1+...+\lambda'_p$
%for any $p$.
%If $T$ and $T'$ are two standard tableaux of shape $Y(\nilp)$ we write $T\preceq T'$
%if $Y_i(T)\preceq Y_i(T')$ for every $i$.

\paragraphe
\label{paragraphe-BuT}
Our construction %of cells decomposition
relies on an $\alpha$-partition of $\B_\nilp$ into subsets parameterized by standard tableaux.
Let us recall the construction, due to Spaltenstein,
of such an $\alpha$-partition.
Let $(V_0,...,V_n)\in\B_\nilp$.
For each $i$ consider the Young diagram $Y(\nilp_{|V_i})$ associated to the restriction $\nilp_{|V_i}$
in the sense of \ref{section-theoreme-principal}.
Let $T$ be a standard tableau of shape $Y(\nilp)$.
The shape of the subtableau of entries $1,...,i$ is a subdiagram
$Y_i(T)\subset Y(\nilp)$ (cf. \ref{standardization}). 
Define $\B_\nilp^T$ as
the set of $\nilp$-stable flags such that $Y(\nilp_{|V_i})=Y_i(T)$ for every $i$.
%the set of $\nilp$-stable flags $(V_0,...,V_n)$ such that 
%the Young diagram $Y(\nilp_{V_i})$ associated to the restriction $\nilp_{|V_i}$
%in the sense of \ref{section-theoreme-principal}
%the Jordan form
%of the restriction $\nilp_{|V_i}$ is represented by the diagram $Y_i(T)$ for every $i$,
%in the sense that the sizes of the Jordan blocks are the lengths of the rows
%of the diagram.
By \cite{Spaltenstein}, the $\B_\nilp^T$'s form an $\alpha$-partition of $\B_\nilp$ into
irreducible, nonsingular subsets of same dimension as $\B_\nilp$.
Therefore,
the components of $\B_\nilp$ are exactly
the closures of the $\B_\nilp^T$'s.
%Here we propose a more general construction.

We generalize this construction. 
Let ${\mathcal R}_n$ denote the set of 
double sequences of integers $(i_k,j_k)_{k=0,...,n}$
with $(i_k)_k$ weakly decreasing, $(j_k)_k$ weakly increasing,
$0\leq i_k\leq j_k\leq n$ and $j_k-i_k=k$ for every $k$.
Let $\rho=(i_k,j_k)_{k}\in{\mathcal R}_n$.
Instead of 
considering the restrictions of $\nilp$ to the subspaces of the flag,
we consider
%arrange the subspaces of the flag into 
the maximal chain of subquotients
\[
0=V_{j_0}/V_{i_0}\subset V_{j_1}/V_{i_1}\subset ...\subset V_{j_n}/V_{i_n}=V
\]
For any $k$ we consider the Young diagram $Y(\nilp_{|V_{j_k}/V_{i_k}})$ associated
to the nilpotent endomorphism of the subquotient $V_{j_k}/V_{i_k}$ induced by $\nilp$.
We denote by $\B_{\nilp,T}^\rho$
the set of $\nilp$-stable flags $(V_0,...,V_n)$ such that %the Jordan form
%of the morphism induced by $\nilp$ on the subquotient $V_{j_k}/V_{i_k}$
%is represented by the diagram $Y_k(T)$.
$Y(\nilp_{|V_{j_k}/V_{i_k}})=Y_k(T)$.
The double sequence $\rho$ being fixed,
we prove that the $\B_{\nilp,T}^\rho$'s 
form an $\alpha$-partition of $\B_\nilp$
(see \ref{BuTrho-alphapartition}) into irreducible, nonsingular
subsets of same dimension as $\B_\nilp$
(see Theorem \ref{theoreme-Spaltenstein-general}).

%Then, %it remains to 
%it is sufficient to
%construct a cell decomposition
%of each $\B_{\nilp,T}^\rho$.
For each $T$,
we construct a cell decomposition
$\B_{\nilp,T}^\rho=\bigcup C^\rho(\tau)$
indexed on row-standard tableaux with $\mathrm{st}(\tau)=T$,
and such that the codimension of $C_\tau^\rho$ in $\B_{\nilp,T}^\rho$
is $n_\mathrm{inv}(\tau)$
(see Theorem \ref{theoreme-decomposition-cellulaire}).

Finally, by collecting together
the cell decompositions of the $\B_{\nilp,T}^\rho$'s for $T$ running over
the set of standard tableaux of shape $Y(\nilp)$ (and fixing $\rho$),
 we get a cell decomposition $\B_\nilp=\bigcup_\tau C^\rho(\tau)$. It 
is not unique, since it depends on the parameter $\rho$.

\paragraphe
Observe that in the cell decomposition of $\B_{\nilp,T}^\rho$ mentioned above,
the dimension of the cells does not depend on $\rho$.
Therefore the cohomology with compact support of $\B_{\nilp,T}^\rho$ only
depends on $T$ (see \ref{bettiBuTrho}).
If $T_\mathrm{min}$ is the minimal standard tableau of shape $Y(\nilp)$
for the dominance order (see \ref{BuTrho-alphapartition}),
then we prove that $\B_{\nilp,T_\mathrm{min}}^\rho$ is a closed subset
of $\B_\nilp$, thus it is a nonsingular irreducible component of $\B_\nilp$.
Then, the cell decomposition allows to compute its Betti numbers.
When $\rho$ is changing, the subset $\B_{\nilp,T_\mathrm{min}}^\rho\subset \B_\nilp$
is changing too, and we get thus a family of components of $\B_\nilp$
which are all nonsingular and have the same Betti numbers.

\paragraphe
Our article contains four parts.
In part 2, we establish some properties of the inversion number $n_\mathrm{inv}(\tau)$.
In geometric part 3,
we prove the results announced in \ref{paragraphe-BuT}.
Part 3 is independent from part 2 before.
In part 4, we apply results of the two parts before to the calculation of Betti numbers.

Fix some conventional notation.
Let $\N=\{0,1,2,...\}$ be the set of nonnegative integers.
Let $\C$ be the field of complex numbers.
Let $\Q$ be the field of rational numbers. 
%We consider the classical cohomology of sheafs.
%Let $H^*(X,\Q)$ denote the cohomology space with rational coefficients
%and let $H_c^*(X,\Q)$ denote the rational cohomology with compact support.
Let $S_n$ be the group of permutations of $\{1,...,n\}$.
We denote by $\# A$ the number of elements in a finite set $A$.
Other pieces of notation will be introduced in what follows.

\partie{Inversion number of row-standard tableaux}

\paragraphe
First we define an elementary operation on row-standard tableaux.
For $i=1,...,n$ let $D_i$ denote the set of row-standard tableaux
$\tau$ which satisfy the following properties
\begin{enumerate}
\item $i$ is not in the first row of $\tau$. Then let $j$ be the entry in the neighbor box above $i$.
\item If $i$ has an entry $i'$ on its right, then $j<i'$.
If $j$ has an entry $j'$ on its right, then $i<j'$.
\item For every $k$
in the same column as $i,j$ and such that $\mathrm{min}(i,j)<k<\mathrm{max}(i,j)$,
either 
$(\mathrm{min}(i,j),k)$ or 
$(k,\mathrm{max}(i,j))$ is an inversion
(but not both).
\end{enumerate}
Let $\tau\in D_i$.
Let $i_1\leq ...\leq i_q=i$ be the entries until $i$ of the row containing $i$,
and let $j_1\leq ...\leq j_q=j$ be the entries until $j$ of the row containing $j$.
Then define $\delta_i(\tau)$ as the tableau obtained by switching $i_k$ and $j_k$
for every $k=1,...,q$. The tableau $\delta_i(\tau)$ remains row-standard.
\ \ Observe that $\delta_i(\tau)\in D_j$ and that we have $\tau=\delta_j\delta_i(\tau)$.
%The map $\delta_i$ so-defined is clearly injective.
%For $i=1,...,n$ we define $D'_i=\{\delta_i(\tau):\tau\in D_i\}$.
%Let $\delta'_i:D'_i\rightarrow D_i$ be the inverse map.

\paragraph{Lemma} \label{elementary-operation} {\em Let $\tau\in D_i$ and let $j$ be the neighbor entry above $i$.
Then we have $n_\mathrm{inv}(\delta_i(\tau))=n_\mathrm{inv}(\tau)-1$
if $(\mathrm{min}(i,j),\mathrm{max}(i,j))$ is an inversion of $\tau$,
and $n_\mathrm{inv}(\delta_i(\tau))=n_\mathrm{inv}(\tau)+1$ otherwise.
%depending on whether $(\mathrm{min}(i,j),\mathrm{max}(i,j))$ is an inversion
%of $\tau$ or not.
}

\medskip
\noindent
{\em Proof.} Let $\mathrm{Inv}(\tau)$ 
(resp. $\mathrm{Inv}(\delta_i(\tau))$)
be the set of inversions of $\tau$ (resp. of $\delta_i(\tau)$).
For $\{k,l\}\cap\{i,j\}=\emptyset$,
it is clear that $(k,l)\in\mathrm{Inv}(\tau)\Leftrightarrow(k,l)\in\mathrm{Inv}(\delta_i(\tau))$.
Let $k$ be in the same column as $i,j$.
Observe that $i,j$ are neighbor in $\tau$ as in $\delta_i(\tau)$, hence 
$k$ is above $i$ if and only if it is above $j$,
and $(k,i)$ have the same relative position in $\tau$ and $\delta_i(\tau)$.
If $k<\mathrm{min}(i,j)$ then it is clear that 
$(k,i)\in\mathrm{Inv}(\tau)\Leftrightarrow(k,j)\in\mathrm{Inv}(\delta_i(\tau))$
and
$(k,j)\in\mathrm{Inv}(\tau)\Leftrightarrow(k,i)\in\mathrm{Inv}(\delta_i(\tau))$.
Likewise if $k>\mathrm{max}(i,j)$ then we have
$(i,k)\in\mathrm{Inv}(\tau)\Leftrightarrow(j,k)\in\mathrm{Inv}(\delta_i(\tau))$
and
$(j,k)\in\mathrm{Inv}(\tau)\Leftrightarrow(i,k)\in\mathrm{Inv}(\delta_i(\tau))$.
Now suppose $\mathrm{min}(i,j)<k<\mathrm{max}(i,j)$.
Say $i<j$ (the other case is treated similarly).
It follows from the definition of inversion that $(i,k)$ is an inversion of $\tau$
if and only if $(k,j)$ is not an inversion of $\delta_i(\tau)$.
Likewise $(k,j)$ is an inversion of $\tau$
if and only if $(i,k)$ is not an inversion of $\delta_i(\tau)$.
By applying condition 3 above, we get
$(i,k)\in\mathrm{Inv}(\tau)\Leftrightarrow (k,j)\notin\mathrm{Inv}(\tau)
\Leftrightarrow (i,k)\in\mathrm{Inv}(\delta_i(\tau))$ and similarly
$(k,j)\in\mathrm{Inv}(\tau)\Leftrightarrow (k,j)\in\mathrm{Inv}(\delta_i(\tau))$. \\
Finally we get that the number of inversions $(k,l)$ with $\{k,l\}\not=\{i,j\}$ is the same for
$\tau$ and $\delta_i(\tau)$. Now observe that, 
as the right-neighbors of $i$ and $j$ are switched 
from $\tau$ to $\delta_i(\tau)$, we have
$(i,j)\in\mathrm{Inv}(\tau)\Leftrightarrow (i,j)\notin \mathrm{Inv}(\delta_i(\tau))$
(resp. $(j,i)\in\mathrm{Inv}(\tau)\Leftrightarrow (j,i)\notin \mathrm{Inv}(\delta_i(\tau))$
if $j<i$). The lemma follows.
\cqfd

\paragraphe Next we show that $n_{\mathrm{inv}}(\tau)$
is the minimal number of operations to transform $\tau$ into its standardization $\mathrm{st}(\tau)$.
We need the following

\paragraph{Lemma} {\em Suppose $\tau\not=\mathrm{st}(\tau)$.
Let $m$ maximal not at the same place in $\tau$ and $\mathrm{st}(\tau)$.
Then $m$ has a below-neighbor entry $i$, which satisfies $i<m$, and we have $\tau\in D_i$ and 
$n_{\mathrm{inv}}(\delta_i(\tau))=n_{\mathrm{inv}}(\tau)-1$.}

\medskip
\noindent{\em Proof.}
By maximality of $m$, %the entries $j>m$ which are in the same column of $m$ 
%are in the last boxes of the column, and 
there is $i'<m$ below $m$, and all $j>m$ of the same column as $m$ are below $i'$.
In particular $m$ has a below-neighbor $i$ and we have $i<m$.
If $m$ has a right neighbor $m'$, then we have $i<m<m'$.
If $i$ has a right neighbor $i'$, then $m'$ also exists,
and by maximality of $m$, the entries in the column of $m'$ 
are in good order from $m'$ to the bottom, in particular
we have $m<m'<i'$.
For $k=i+1,...,m-1$ in the same column as $i,m$, 
either $k$ is above $i,m$, then $(i,k)$ is an inversion and $(k,m)$ is not one,
or $k$ is below $i,m$, then $(i,k)$ is not an inversion and $(k,m)$ is one.
Therefore $\tau\in D_i$ and $(i,m)$ is an inversion.
By Lemma \ref{elementary-operation} it follows $n_\mathrm{inv}(\delta_i(\tau))=n_\mathrm{inv}(\tau)-1$.
\cqfd

\paragraph{Proposition} \label{proposition-coxeter} {\em Let $\tau$ be row-standard, then there is a sequence of integers
$i_1,...,i_m$ such that $\mathrm{st}(\tau)=\delta_{i_1}\cdots\delta_{i_m}(\tau)$.
The inversion number $n_{\mathrm{inv}(\tau)}$ is the minimal length
of such a sequence.}

\medskip
\noindent
{\em Proof.} 
If there are $i_1,...,i_m$ with $\mathrm{st}(\tau)=\delta_{i_1}\cdots\delta_{i_m}(\tau)$,
then we get $m\geq \mathrm{inv}(\tau)$ by Lemma \ref{elementary-operation}.
We prove by induction on $n_{\mathrm{inv}}(\tau)$ 
that there are $i_1,...,i_m$ with $m=n_{\mathrm{inv}}(\tau)$
such that
$\mathrm{st}(\tau)=\delta_{i_1}\cdots\delta_{i_m}(\tau)$.
If $n_{\mathrm{inv}}(\tau)=0$ then $\tau=\mathrm{st}(\tau)$
and it is immediate.
Suppose $n_{\mathrm{inv}}(\tau)>0$. 
By the lemma above there is $i$ such that $\tau\in D_i$ and 
$n_{\mathrm{inv}}(\delta_i(\tau))=n_{\mathrm{inv}}(\tau)-1$.
The property follows by induction hypothesis applied to $\delta_i(\tau)$. 
\cqfd

\medskip

We construct a graph whose vertices are row-standard tableaux of shape $Y(\nilp)$
and with one edge between $\tau$ and $\tau'$ if there is $i$ such that $\tau'=\delta_i(\tau)$.\ \
For $Y(\nilp)=\mbox{\scriptsize$\yng(2,2,1)$}$, we get for example the following graph.

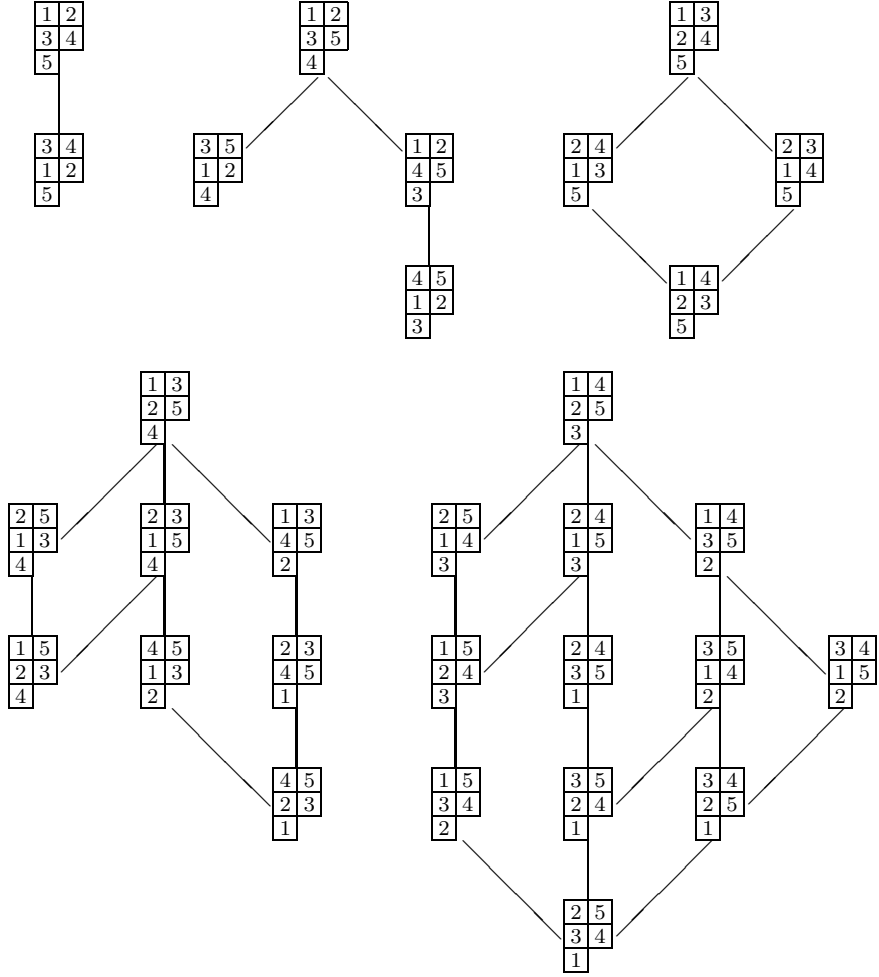
\begin{picture}(350,380)(0,10)

\put(40,370){\mbox{\scriptsize $\young(12,34,5)$}}
\put(40,320){\mbox{\scriptsize $\young(34,12,5)$}}
\put(49,358){\line(0,-1){22}}

\put(140,370){\mbox{\scriptsize $\young(12,35,4)$}}
\put(100,320){\mbox{\scriptsize $\young(35,12,4)$}}
\put(180,320){\mbox{\scriptsize $\young(12,45,3)$}}
\put(151,357){\line(1,-1){28}}
\put(147,357){\line(-1,-1){27}}
\put(180,270){\mbox{\scriptsize $\young(45,12,3)$}}
\put(189,308){\line(0,-1){22}}

\put(280,370){\mbox{\scriptsize $\young(13,24,5)$}}
\put(240,320){\mbox{\scriptsize $\young(24,13,5)$}}
\put(320,320){\mbox{\scriptsize $\young(23,14,5)$}}
\put(291,357){\line(1,-1){28}}
\put(287,357){\line(-1,-1){27}}
\put(280,270){\mbox{\scriptsize $\young(14,23,5)$}}
\put(251,307){\line(1,-1){28}}
\put(327,307){\line(-1,-1){27}}

\put(240,230){\mbox{\scriptsize $\young(14,25,3)$}}
\put(190,180){\mbox{\scriptsize $\young(25,14,3)$}}
\put(240,180){\mbox{\scriptsize $\young(24,15,3)$}}
\put(290,180){\mbox{\scriptsize $\young(14,35,2)$}}
\put(246,218){\line(-1,-1){36}}
\put(249,218){\line(0,-1){22}}
\put(252,218){\line(1,-1){36}}
\put(190,130){\mbox{\scriptsize $\young(15,24,3)$}}
\put(240,130){\mbox{\scriptsize $\young(24,35,1)$}}
\put(290,130){\mbox{\scriptsize $\young(35,14,2)$}}
\put(340,130){\mbox{\scriptsize $\young(34,15,2)$}}
\put(246,168){\line(-1,-1){36}}
\put(302,168){\line(1,-1){37}}
\put(199,168){\line(0,-1){22}}
\put(249,168){\line(0,-1){22}}
\put(299,168){\line(0,-1){22}}
\put(190,80){\mbox{\scriptsize $\young(15,34,2)$}}
\put(240,80){\mbox{\scriptsize $\young(35,24,1)$}}
\put(290,80){\mbox{\scriptsize $\young(34,25,1)$}}
\put(296,118){\line(-1,-1){36}}
\put(346,118){\line(-1,-1){36}}
\put(199,118){\line(0,-1){22}}
\put(249,118){\line(0,-1){22}}
\put(299,118){\line(0,-1){22}}
\put(240,30){\mbox{\scriptsize $\young(25,34,1)$}}
\put(202,68){\line(1,-1){37}}
\put(296,68){\line(-1,-1){36}}
\put(249,68){\line(0,-1){22}}

\put(80,230){\mbox{\scriptsize $\young(13,25,4)$}}
\put(30,180){\mbox{\scriptsize $\young(25,13,4)$}}
\put(80,180){\mbox{\scriptsize $\young(23,15,4)$}}
\put(130,180){\mbox{\scriptsize $\young(13,45,2)$}}
\put(92,218){\line(1,-1){37}}
\put(86,218){\line(-1,-1){36}}
%\put(91.5,218){\line(1,-1){38}}
\put(89,218){\line(0,-1){22}}
\put(30,130){\mbox{\scriptsize $\young(15,23,4)$}}
\put(80,130){\mbox{\scriptsize $\young(45,13,2)$}}
\put(130,130){\mbox{\scriptsize $\young(23,45,1)$}}
\put(39,168){\line(0,-1){22}}
\put(89,168){\line(0,-1){22}}
\put(139,168){\line(0,-1){22}}
\put(86,168){\line(-1,-1){36}}
\put(130,80){\mbox{\scriptsize $\young(45,23,1)$}}
\put(139,118){\line(0,-1){22}}
\put(92,118){\line(1,-1){37}}

\end{picture}

\noindent
Each connected component contains a unique standard tableau.
Two tableaux $\tau$ and $\tau'$ are in the same connected component
if we have $\mathrm{st}(\tau)=\mathrm{st}(\tau')$.
By the proposition, the number of inversion $n_\mathrm{inv}(\tau)$ is the length
between $\tau$ and $\mathrm{st}(\tau)$ in the graph.

\paragraphe
\label{distribution-longueur}
Let $T$ be standard. For each $i$, let $q_i$ be the number of the column of $T$ containing $i$
and let $p_i$ be the number of rows of length $q_i$ in the subtableau $T[1,...i]$.
The next proposition allows to describe the distribution of inversion numbers.
 %First we prove the following

%\paragraph{Lemma} {\em %(a) 
%If $\tau\in D_i$ and $i+1,...,n$ have the same place in $\tau$ and $\mathrm{st}(\tau)$,
%then $i+1,...,n$ remain at the same place in $\delta_i(\tau)$
%and $\mathrm{st}(\tau)$.}
%If $i$ is not in the first row of $\tau$
%and if $i+1,...,n$ have the same place in $\tau$ and $\mathrm{st}(\tau)$,
%then $\tau\in D_i$ and $n_\mathrm{inv}(\delta_i(\tau))=n_\mathrm{inv}(\tau)+1$. 
%Moreover $i+1,...,n$ remain at the same place in $\delta_i(\tau)$
%and $\mathrm{st}(\tau)$.}% \\
%(b) If $i$ is maximal which has not the same place in $\tau$ and $\mathrm{st}(\tau)$,
%then $\tau\in D'_i$ and $n_\mathrm{inv}(\delta'_i(\tau))=n_\mathrm{inv}(\tau)-1$.}

%\medskip
%\noindent
%{\em Proof.}
%Every $k\in\{i+1,...,n\}$ of the column of $i$ is at the same place in $\tau$ and $\mathrm{st}(\tau)$,
%hence $i$ is above each such $k$. 
%If $j$ denotes the entry just above $i$, then it follows $j<i$. 
%Thus the entries deplaced from $\tau$ to $\delta_i(\tau)$
%are all weakly lower than $i$.
%\cqfd

\paragraph{Proposition} {\em %Let $T$ be standard. Let $p_i$ be the number of rows 
%in the subtableau $T[1,...,i]$ which have the same length as the row of $i$. \\
(a) Let $\kappa_1,...,\kappa_n$ be integers with $0\leq \kappa_i\leq p_i-1$ for any $i$.
Then $(\delta_n)^{\kappa_n}\cdots (\delta_1)^{\kappa_1}(T)$ is well-defined. \\
(b) For every $\tau$ row-standard such that $\mathrm{st}(\tau)=T$, there are unique
integers $\kappa_1,...,\kappa_n$ with $0\leq \kappa_i\leq p_i-1$
such that we have $\tau=(\delta_n)^{\kappa_n}\cdots (\delta_1)^{\kappa_1}(T)$.
Moreover $n_\mathrm{inv}(\tau)=\kappa_1+...+\kappa_n$.}

\medskip
By the proposition, we obtain the formula \\
$
%\begin{array}{l}
\#\{\tau:\mathrm{st}(\tau)=T,\ n_\mathrm{inv}(\tau)=m\} %\\
%\qquad\qquad
=\#\{(\kappa_1,...,\kappa_n):
0\leq \kappa_i\leq p_i-1\,\forall i,\ \sum_{i=1}^n\kappa_i=m\}.
%\end{array}
$ 

\bigskip
\noindent
{\em Proof.} 
First, observe that,
if $\tau\in D_i$ and if $i+1,...,n$ have the same place in $\tau$ and $\mathrm{st}(\tau)$,
then $i+1,...,n$ remain at the same place in $\delta_i(\tau)$
and $\mathrm{st}(\tau)$.
Indeed every $k\in\{i+1,...,n\}$ in the column of $i$ is then
either below or on the right of $i$, hence it keeps the same place
in $\delta_i(\tau)$.

(a) We reason by induction on $n$ with immediate initialization in $1$.
Let us prove the property for $n\geq 2$.
By induction hypothesis (considering the subtableau of entries $1,...,n-1$)
the tableau $\tau'=(\delta_{n-1})^{\kappa_{n-1}}\cdots (\delta_1)^{\kappa_1}(T)$ is well defined and
$n$ has the same place in $\tau'$ and $T$.
Then we reason by induction on $\kappa_n\geq 0$ with immediate initialization for $\kappa_n=0$.
Let us prove the property for $\kappa_n\geq 1$.
By induction hypothesis, $\tau''=(\delta_n)^{\kappa_n-1}(\tau')$ is well defined.
The entry $n$ has been moved by $\kappa_n-1$ ranks to the up from $\tau'$ to $\tau''$.
As $\kappa_n<p_n$, there is an entry $j$ just above $n$ in $\tau''$, and $j$ is the last box of its own row.
As in the proof of Lemma \ref{proposition-coxeter}, each
$k=j+1,...,n-1$ in the same column as $j,n$ is such that either $(j,k)$ or $(k,n)$ is an inversion
(but not both), therefore $\tau''\in D_n$, and $\delta_n(\tau'')$ is well-defined.

(b) %We prove the uniqueness by induction on $n$ with immediate initialization for $n=1$.
%Let us prove the property for $n\geq 2$.
Suppose 
$\tau=(\delta_n)^{\kappa_n}\cdots (\delta_1)^{\kappa_1}(T)
=(\delta_n)^{\kappa'_n}\cdots (\delta_1)^{\kappa'_1}(T)$.
Then $\kappa_n$ (and similarly $\kappa'_n$) is the number of boxes 
below $n$ in $\tau$.
Thus $\kappa_n=\kappa'_n$.
As $\delta_n$ is injective we get $(\delta_{n-1})^{\kappa_{n-1}}\cdots (\delta_1)^{\kappa_1}(T)
=(\delta_{n-1})^{\kappa'_{n-1}}\cdots (\delta_1)^{\kappa'_1}(T)$.
Then $\kappa_{n-1}$ (and similarly $\kappa'_{n-1}$) is the number of entries $j<n$ below
$n-1$ in this new tableau. Thus $\kappa_{n-1}=\kappa'_{n-1}$.
And so on...
We deduce $\kappa_i=\kappa'_i$ for any $i$.
Moreover we see that, if $m+1,...,n$ have the same place in $\tau$ and $T$, then we must
have $\kappa_n=...=\kappa_{m+1}=0$. \\
We prove the existence by induction on $\mathrm{inv}(\tau)$ with immediate initialization for 
$\mathrm{inv}(\tau)=0$.
Suppose $\mathrm{inv}(\tau)>0$. Then there is $m$ maximal which has not the same
place in $\tau$ and $T$. By Lemma \ref{proposition-coxeter} there is $i<m$ just below $m$
and we have $\tau\in D_i$ and $n_{\mathrm{inv}}(\delta_i(\tau))=n_\mathrm{inv}(\tau)-1$.
Let $\tau'=\delta_i(\tau)$. Then $\tau'\in D_m$ and $\tau=\delta_m(\tau')$.
By induction hypothesis
%and as $m+1,...,n$ have the same place in $\tau'$ and $T$,
we have $\tau'=(\delta_m)^{\kappa'_m}\cdots (\delta_1)^{\kappa'_1}(T)$
with $\kappa'_1+...+\kappa'_m=n_\mathrm{inv}(\tau')$.
We get $\tau=(\delta_m)^{\kappa'_m+1}\cdots (\delta_1)^{\kappa'_1}(T)$
and we have $\kappa'_1+...+(\kappa'_m+1)=n_\mathrm{inv}(\tau)$.
\cqfd

\partie{Geometric constructions}

\paragraphe \label{BuTrho-alphapartition}
We deal with the partition $\B_\nilp=\bigsqcup_T\B_{\nilp,T}^\rho$ introduced in \ref{paragraphe-BuT}.
Let us recall the dominance order on standard tableaux.
For $T$ standard,
let $c_{\leq q}T[1,...,i]$
be the number of boxes in the first $q$ columns of the subtableau of entries $1,...,i$.
We write $T\preceq T'$ if $c_{\leq q}T[1,...,i]\geq c_{\leq q}T'[1,...,i]$ for any $i$ and any $q$. First we prove the following

\paragraph{Proposition}{\em Fix $\rho\in{\mathcal R}_n$.
Let $T$ be standard. Then we have $\overline{\B_{\nilp,T}^\rho}\subset \bigcup_{S\preceq T}\B_{\nilp,S}^\rho$
where the union is taken over standard tableaux $S$ such that $S\preceq T$.}

\bigskip

It follows from this proposition that the $\B_{\nilp,T}^\rho$ for $\rho$ fixed and $T$ running over the set
of standard tableaux of shape $Y(\nilp)$ form an $\alpha$-partition of $\B_\nilp$.
Indeed, take a total order on standard tableaux completing the dominance order.
Then, the $\B_{\nilp,T}^\rho$'s, arranged according to this order, form a sequence whose first $k$ terms
always have a closed union.

There is a unique tableau $T_\mathrm{min}$ of shape $Y(\nilp)$ which is minimal for the dominance order.
Let $\mu_1,...,\mu_s$ be the lengths of the columns of $Y(\nilp)$.
Then $T_\mathrm{min}$ is the standard tableau with numbers $1,...,\mu_1$ in the first column,
$\mu_1+1,...,\mu_1+\mu_2$ in the second column, etc.
For example for $\lambda(\nilp)=(3,2,2)$ we have
\[T_\mathrm{min}=\young(147,25,36)\]
By the proposition we get that $\B_{\nilp,T_\mathrm{min}}^\rho$ is a closed subset of $\B_\nilp$.

\medskip
\noindent
{\em Proof of the proposition.}
By definition, for $(V_0,...,V_n)\in\B_{\nilp,T}^\rho$
and $(i_k,j_k)\in\rho$ the %sizes of the Jordan blocks of $\nilp_{|V_{j_k}/V_{i_k}}$
Young diagram $Y(\nilp_{|V_{j_k}/V_{i_k}})$ associated to the nilpotent map induced by $\nilp$ on the subquotient 
$V_{j_k}/V_{i_k}$
is the shape 
%are the lengths of the rows 
of the subtableau $T[1,...,k]$.
Thus the number of boxes in the first $q$ columns of both
coincide for any $q\geq 1$.
%Let $c_{\leq q}T[1,...,k]$ be the number of boxes in the first $q$ columns
%of the subtableau $T[1,...,k]$. 
Thus $\mathrm{dim}\,\mathrm{ker}\,\nilp^q_{|V_{j_k}/V_{i_k}}=c_{\leq q}T[1,...,k]$. \\
Suppose $\B_{\nilp,S}^\rho\cap \overline{\B_{\nilp,T}^\rho}$ nonempty and take
$(V_0,...,V_n)\in\B_{\nilp,S}^\rho\cap \overline{\B_{\nilp,T}^\rho}$.
Then we have $\mathrm{dim}\,\mathrm{ker}\,\nilp^q_{|V_{j_k}/V_{i_k}}\geq c_{\leq q}T[1,...,k]$ by the following lemma.
It follows $c_{\leq q}S[1,...,k]\geq c_{\leq q}T[1,...,k]$ for any $q\geq 1$ and $k=1,...,n$,
therefore $S\preceq T$.
\cqfd

\paragraph{Lemma}{\em 
The set $\{(V_0,...,V_n)\in \B_\nilp:\mathrm{dim}\,\mathrm{ker}\,\nilp^q_{|V_j/V_i}\geq c\}$
is closed for any $c\in\N$.}

\medskip
\noindent
{\em Proof.}
We show the lemma for $q=1$. For the general case, replace $\nilp$ by $\nilp^q$.
We prove that $\mathrm{dim}\,\mathrm{ker}\,\nilp_{|V_j/V_i}=j+i-\mathrm{dim}\,(V_i+\nilp(V_j))$.
Then the lemma follows from the lower semicontinuity of the map
$(V_i,V_j)\mapsto \mathrm{dim}\,(V_i+\nilp(V_j))$ defined on the product of grassmannians. \\
We have $\mathrm{ker}\,\nilp_{|V_{j}/V_{i}}=\nilp^{-1}(V_i)\cap V_j$.
By the rank formula applied to the restriction of $\nilp$ to
$\nilp^{-1}(V_i)\cap V_j$ 
we get $\mathrm{dim}\,\nilp^{-1}(V_i)\cap V_j=\mathrm{dim}\,V_j\cap\mathrm{ker}\,\nilp+\mathrm{dim}\,V_i\cap \nilp(V_j)$.
On one hand, we have $\mathrm{dim}\,V_i\cap
\nilp(V_j)=i+\mathrm{dim}\,\nilp(V_j)-\mathrm{dim}\,(V_i+\nilp(V_j))$. On the
other hand, the rank formula gives
$\mathrm{dim}\,V_j\cap\mathrm{ker}\,\nilp=j-\mathrm{dim}\,\nilp(V_j)$. The desired formula follows.
\cqfd

\paragraphe
\label{theoreme-Spaltenstein-general}
The following theorem
generalizes \cite[\S II.5.5]{Spaltenstein}.

\paragraph{Theorem} {\em Fix $\rho\in {\mathcal R}_n$. Let $T$ be standard.
The set $\B_{\nilp,T}^\rho$ is an irreducible, nonsingular
subvariety of $\B_\nilp$ and we have $\mathrm{dim}\,\B_{\nilp,T}^\rho=\mathrm{dim}\,\B_\nilp$.}

\bigskip

The theorem is proved by induction
in sections \ref{section-preuve-1}--\ref{section-preuve-5}. 
From the theorem and Proposition \ref{BuTrho-alphapartition}, we deduce the following

\paragraph{Corollary} {\em Fix $\rho\in{\mathcal R}_n$. For every $T$, 
the closure $\overline{\B_{\nilp,T}^\rho}$ is an irreducible component
of $\B_\nilp$ and every irreducible component is obtained in that way.
Moreover we have $\overline{\B_{\nilp,T_\mathrm{min}}^\rho}=\B_{\nilp,T_\mathrm{min}}^\rho$ and
in particular this component is nonsingular.}

\bigskip

For each $\rho$, we obtain a different parameterization of the components
of $\B_\nilp$ by standard tableaux, and the $\B_{\nilp,T_\mathrm{min}}^\rho$'s
for $\rho$ running over the set ${\mathcal R}_n$ form a family of nonsingular components.

\medskip
\noindent
{\em Remark.}
Let us describe the link between the different parameterizations of the components.
Take as reference the component ${\mathcal K}^T=\overline{\B_\nilp^T}$
obtained as the closure of the Spaltenstein set (see \ref{paragraphe-BuT}),
and let us describe $S$ such that ${\mathcal K}^T=\overline{\B_{\nilp,S}^\rho}$.
By \cite[Theorem 3.3]{vanLeeuwen}, 
for $(V_0,...,V_n)\in {\mathcal K}^T$ generic and any $1\leq i<j\leq n$
the Young diagram 
$Y(\nilp_{|V_j/V_i})$ is the shape of the
tableau obtained as rectification by jeu de taquin
of the subtableau $T[i+1,...,j]$.
%(see \cite{vanLeeuwen} for details about jeu de taquin).
Write $\rho=(i_k,j_k)_k$.
For each $k$ let $Y^{(k)}$ be the Young diagram forming the shape
of the rectification by jeu de taquin
of the subtableau $T[i_k+1,...,j_k]$.
We get a chain of subdiagrams $\emptyset=Y^{(0)}\subset Y^{(1)}\subset  ...\subset Y^{(n)}=Y(\nilp)$
and $S$ is the standard tableau which corresponds to this chain
in the sense of \ref{standardization}.

\paragraphe
\label{theoreme-decomposition-cellulaire}
We fix $\rho\in{\mathcal R}_n$.
The main result of this section states the existence of a cell decomposition for each $\B_{\nilp,T}^\rho$.

\paragraph{Theorem} {\em
Let $T$ be standard. The set $\B_{\nilp,T}^\rho$ has a cell decomposition 
$\B_{\nilp,T}^\rho=\bigsqcup C^\rho(\tau)$ parameterized
by the row-standard tableaux $\tau$ of standardization $\mathrm{st}(\tau)=T$,
such that the codimension of the cell $C^\rho(\tau)$ in $\B_{\nilp,T}^\rho$
is equal to the inversion number $n_{\mathrm{inv}}(\tau)$.}

\bigskip

As said in \ref{BuTrho-alphapartition} the subsets $\B_{\nilp,T}^\rho$ form
an $\alpha$-partition of $\B_\nilp$. Therefore by collecting together
the cell decompositions of the $\B_{\nilp,T}^\rho$'s for $T$ running over the set of standard tableaux,
we obtain a cell decomposition of $\B_\nilp$. This proves Theorem \ref{section-theoreme-principal}.

\bigskip

We prove both theorems simultaneously, by induction on the dimension of $V$.

\subsubsection*{Proof of Theorems \ref{theoreme-Spaltenstein-general}
and \ref{theoreme-decomposition-cellulaire}}

\paragraphe
\label{section-preuve-1}
First, we point out a duality
in the family parameterized by $\rho\in{\mathcal R}_n$ of partitions of $\B_\nilp$.
It will allow us to suppose that the sequence $\rho=(i_k,j_k)_k$ is such that $(i_{n-1},j_{n-1})=(0,n-1)$.

Let $V^*$ be the dual vector space of $V$.
The dual map $\nilp^*:V^*\rightarrow V^*$ is also nilpotent.
Let $\B_{\nilp^*}$ be the Springer fiber relative to $\nilp^*$.
The maps $\nilp^*$ and $\nilp$ are conjugated, %under some isomorphism $g:V\rightarrow V^*$,
in particular they have the same Jordan form.
%The map
%$\B_\nilp\rightarrow \B_{\nilp^*},\ (V_0,...,V_n)\mapsto (gV_n,...,gV_0)$
%is an isomorphism of algebraic varieties, which restricts to an isomorphism
%between the subsets $\B_{\nilp,T}^\rho$ and $\B_{\nilp^*,T}^\rho$ for any $\rho\in{\mathcal R}_n$
%and any standard tableau $T$ of shape $Y(\nilp)$.
For a subspace $W\subset V$ let $W^\perp=\{\phi\in V^*:\phi(w)=0\ \forall w\in W\}$.
The map
\[
\Psi:\B_\nilp\rightarrow \B_{\nilp^*},\ (V_0,...,V_n)\mapsto (V_n^\perp,...,V_0^\perp)\]
is well-defined and is an isomorphism of algebraic varieties.
Writing $\rho=(i_k,j_k)_{k=0,...,n}$,
let us define $\rho^*=\left(i^*_k,j^*_k\right)_{k=0,...,n}\in{\mathcal R}_n$
by
$i^*_k=n-j_k$ and $j^*_k=n-i_k$ for every $k$.
%
%\paragraph{Lemma} {\em
The map $\Psi$ restricts to an isomorphism of algebraic varieties
between $\B_{\nilp,T}^\rho$ and $\B_{\nilp^*,T}^{\rho^*}$,
for every standard tableau $T$. 
Indeed,
for $\F=(V_0,...,V_n)\in\B_\nilp$ and
any $k=1,...,n$, the quotient
$V_{i_k}^\perp/V_{j_k}^\perp$ is naturally isomorphic to the dual space $(V_{j_k}/V_{i_k})^*$,
and the endomorphism
$(\nilp^*)_{|V_{i_k}^\perp/V_{j_k}^\perp}$ induced by $\nilp^*$ coincides with the dual map of
$\nilp_{|V_{j_k}/V_{i_k}}$.
It follows that the linear maps $(\nilp^*)_{|V_{i_k}^\perp/V_{j_k}^\perp}$ and $\nilp_{|V_{j_k}/V_{i_k}}$
are conjugated, thus they have the same Jordan form.
Therefore, we have $\Psi(\B_{\nilp,T}^\rho)=\B_{\nilp^*,T}^{\rho^*}$
for every $T$.

%By
%combining these observations,
%we get that the subvarieties $\B_{\nilp,T}^\rho$
%and $\B_{\nilp,T}^{\rho^*}$ are isomorphic
%for every $T$.
In what follows, we may thus assume that $\rho=(i_k,j_k)_k$ is such that
$(i_{n-1},j_{n-1})=(0,n-1)$,
since otherwise we can %consider $\rho^*$
deal with $(\nilp^*,\rho^*)$ instead of $(\nilp,\rho)$.

\paragraphe
\label{section-preuve-2}
Let 
${\mathcal H}_\nilp$ be the set of $\nilp$-stable hyperplanes $H\subset V$.
Let $Z(\nilp)\subset GL(V)$ be the (closed) subgroup of elements
$g$ such that $g\nilp=\nilp g$.
The group $Z(\nilp)$ is connected since it is an open subset
of the vector space of endomorphisms which commute with $\nilp$.
The action of $Z(\nilp)$ on hyperplanes leaves ${\mathcal H}_\nilp$ invariant.
The action of $Z(\nilp)$ on flags leaves the Springer fiber $\B_\nilp$ invariant.
The map 
\[\Phi:\B_{\nilp}\rightarrow {\mathcal H}_\nilp,\ \ (V_0,...,V_n)\mapsto V_{n-1}\]
is algebraic and $Z(\nilp)$-equivariant.

Now we fix a standard tableau $T$.
It is easy to see that the action of $Z(\nilp)$ on flags leaves $\B_{\nilp,T}^\rho$ invariant.
We consider the restriction of $\Phi$ to $\B_{\nilp,T}^\rho$
\[\Phi_T:\B_{\nilp,T}^\rho\rightarrow {\mathcal H}_\nilp,\ \ (V_0,...,V_n)\mapsto V_{n-1}\]
which is algebraic and $Z(\nilp)$-equivariant.
Let $T'$ be the subtableau obtained from $T$ by deleting the box number $n$.
Let $Y'$ be the shape of $T'$, which is the subdiagram of $Y(\nilp)$
obtained by deleting the same box.
Write $\rho'=(i_k,j_k)_{k=1,...,n-1}$.
The image of $\Phi_T$ is the subset of $\nilp$-stable hyperplanes $H$
such that the Young diagram $Y(\nilp_{|H})$ associated to the 
restriction of $\nilp$ to $H$ is equal to $Y'$.
Let $H\in{\mathcal H}_\nilp$ be such a hyperplane.
Then we have
\[\Phi_T^{-1}(H)=\{(V_0,...,V_n)\in\B_{\nilp,T}^\rho:V_{n-1}=H\}=\B_{\nilp_{|H},T'}^{\rho'}\]
where $\B_{\nilp_{|H},T'}^{\rho'}$ is the subset which corresponds to $T'$ in the Springer fiber
$\B_{\nilp_{|H}}$ associated to the nilpotent map $\nilp_{|H}:H\rightarrow H$.

%The map
%\[Z(\nilp)\times \Phi^{-1}(H)\rightarrow \B_{\nilp,T}^\rho,\ \ 
%g,(V_0,...,V_{n-1})\mapsto (gV_0,...,gV_{n-1},V_n)\]
%is well-defined and algebraic.
We prove Theorems \ref{theoreme-Spaltenstein-general} and \ref{theoreme-decomposition-cellulaire}
by induction on $n=\mathrm{dim}\,V$. For Theorem \ref{theoreme-Spaltenstein-general} we
show that $\Phi_T$ is locally trivial.
For Theorem \ref{theoreme-decomposition-cellulaire},
using the local triviality of $\Phi_T$,
we construct a cell decomposition of $\B_{\nilp,T}^\rho$
over a cell decomposition of the image of $\Phi_T$.

\paragraphe
\label{section-preuve-3}
First, we study the action of $Z(\nilp)$ on ${\mathcal H}_\nilp$.
Note that a hyperplane $H$ is $\nilp$-stable if and only if $H\supset \mathrm{Im}\,\nilp$.
Let $W=V/\mathrm{Im}\,\nilp$ and let $\zeta:V\rightarrow W$ be the
surjective linear map.
Then the variety ${\mathcal H}_\nilp$ is isomorphic to the variety
${\mathcal H}(W)$ of hyperplanes of $W$.
Each $g\in Z(\nilp)$ defines a quotient map in $GL(W)$.
We get a morphism of algebraic groups $\varphi:Z(\nilp)\rightarrow GL(W)$.
Then $Z(\nilp)$ acts linearly on ${\mathcal H}(W)$ and the isomorphism
${\mathcal H}_\nilp\cong {\mathcal H}(W)$ is $Z(\nilp)$-equivariant.

The iterated kernels form a partial flag 
$0\subset \mathrm{ker}\,\nilp\subset \mathrm{ker}\,\nilp^2\subset ...\subset \mathrm{ker}\,\nilp^s=V$.
We get a partial flag of $W$:
\[0\subset \zeta(\mathrm{ker}\,\nilp)\subset \zeta(\mathrm{ker}\,\nilp^2)\subset ...\subset \zeta(\mathrm{ker}\,\nilp^s)=W.\]
Let $W_q=\zeta(\mathrm{ker}\,\nilp^q)$.
Let 
\[P=\{g\in GL(W):g(W_q)=W_q\ \forall q\}.\]
This is a parabolic subgroup.
It is easy to see that each kernel $\mathrm{ker}\,\nilp^q$ is invariant by $g\in Z(\nilp)$,
hence the image of $\varphi$ is contained in $P$.
We prove the following

\paragraph{Lemma}{\em There is a morphism of algebraic groups
$\psi:P\rightarrow Z(\nilp)$ such that $\varphi\circ\psi=\mathrm{id}_P$.}

\medskip
\noindent
{\em Proof.}
We fix a linear embedding $\xi:W\hookrightarrow V$ such that $\zeta\circ \xi=\mathrm{id}_{W}$
and such that in addition $\xi(W_q)\subset \mathrm{ker}\,\nilp^q$.
Hence $\xi(W_q)=\xi(W)\cap \mathrm{ker}\,\nilp^q$.
Any $g\in P$ induces a linear map $\xi g\xi^{-1}:\xi(W)\rightarrow V$.
Let $W'=\xi(W)$ and $g'=\xi g\xi^{-1}$.
Let us prove that there is a unique linear map
$\overline{g}:V\rightarrow V$ commuting with $\nilp$
which extends $g'$.
We have $V=\bigoplus_{q=1}^{s-1}\nilp^q(W')$.
For $v=\nilp^q(w)\in \nilp^q(W')$, we must have $\overline{g}(v)=\nilp^q(g'(w))$.
Thus the extension is unique.
Let us show that $\overline{g}$ defined in that way on $\nilp^q(W')$ is well-defined.
If $v=\nilp^q(w)=\nilp^q(w')$ with $w,w'\in W'$ then $w-w'\in\mathrm{ker}\,\nilp^q$.
As $g$ leaves $W_q$ invariant and as $\xi(W_q)=W'\cap\mathrm{ker}\,\nilp^q$,
we get $g'(w-w')\in\mathrm{ker}\,\nilp^q$ hence $\nilp^qg'(w-w')=0$.
Thus $\overline{g}(v)=\nilp^q(g'(w))=\nilp^q(g'(w'))$ is well-defined.
It is straightforward to show that the map so obtained is linear on $\nilp^q(W')$
and it follows from the definition that $\overline{g}\nilp=\nilp\overline{g}$ on $\nilp^q(W')$.
By collecting together these maps on the $\nilp^q(W')$'s we get a linear map $\overline{g}:V\rightarrow V$
which commutes with $\nilp$. 
By construction, the map $g\mapsto\overline{g}$ is algebraic.
Moreover, by uniqueness, we have 
$\overline{g}\circ \overline{g^{-1}}=I$
and $\overline{h\circ g}=\overline{h}\circ\overline{g}$ for $g,h\in P$.
Therefore, the map $\psi:P\rightarrow Z(\nilp)$ defined by $\psi(g)=\overline{g}$
is a morphism of algebraic groups.
\cqfd

\bigskip
By the lemma
the orbits of ${\mathcal H}(W)$ for the action of $Z(\nilp)$ are the orbits for
the action of $P$, which are the subsets
${\mathcal H}(W)_q$ defined for $q=1,...,s$ by 
\[{\mathcal H}(W)_q=\{H\in {\mathcal H}(W): H\supset W_{q-1}\mbox{ and }H\not\supset W_{q}\}.\]
The orbits of ${\mathcal H}_\nilp$ for the action of $Z(\nilp)$ are thus the subsets
${\mathcal H}_{\nilp,q}$ defined for $q=1,...,s$ by 
\[{\mathcal H}_{\nilp,q}=\{H\in {\mathcal H}_\nilp: H\supset \mathrm{ker}\,\nilp^{q-1}\mbox{ and }H\not\supset 
\mathrm{ker}\,\nilp^{q}\}.\]
Recall that $\lambda(\nilp)=(\lambda_1\geq ...\geq \lambda_r)$ 
is the partition of $n$ formed by
the sizes of the
Jordan blocks of $\nilp$, and $Y(\nilp)$ is the Young diagram of rows of lengths
$\lambda_1,...,\lambda_r$.
Let $\mu=(\mu_1\geq ...\geq\mu_s)$ be the conjugated partition.
Thus $\mu_q$ is the length of the $q$-th column of the diagram.
In particular $\mu_1=r$ is the length of the first row.
Let $\mu'_q=\mu_q-\mu_{q+1}$ for $q<s$ and $\mu'_s=\mu_s$. This is the number of rows of length $q$ in the diagram.
We have $\mathrm{dim}\,W_q=\mu'_1+...+\mu'_q$.
Observe that ${\mathcal H}(W)_{q}$ is isomorphic to an open subset of ${\mathcal H}(W/W_{q-1})$,
the variety of hyperplanes of $W/W_{q-1}$.
Therefore $\mathrm{dim}\,{\mathcal H}(W)_{q}=\mu_q-1$.

Let $B\subset P$ be a Borel subgroup.
%It and all its subgroups can be seen as subgroups of $Z(\nilp)$.
The orbits of ${\mathcal H}_\nilp$ for the action of $B$
form a cell decomposition of ${\mathcal H}_\nilp$,
which can be written ${\mathcal H}_\nilp=\bigsqcup_{l=1}^r{\mathcal C}(l)$
so that $\mathrm{dim}\,{\mathcal C}(l)=l-1$
and ${\mathcal H}_{\nilp,q}={\mathcal C}(\mu_{q+1}+1)\sqcup...\sqcup{\mathcal C}(\mu_{q})$.
For each $l$ choose a representative $H_l\in {\mathcal C}(l)$.
There is a unipotent subgroup $U(l)\subset B$
such that the map $\phi_l:U(l)\rightarrow C(l),\ g\mapsto gH_l$
is an isomorphism of algebraic varieties.
(We use the isomorphism ${\mathcal H}_\nilp\cong {\mathcal H}(W)$
and we consider the Schubert cell decomposition of ${\mathcal H}(W)$, see \cite[\S 1.1]{Brion}).
In particular ${\mathcal C}(\mu_{q})$ is an open subset of the orbit ${\mathcal H}_{\nilp,q}$.

\paragraphe
\label{section-preuve-4}
We come back to the map $\Phi_T:\B_{\nilp,T}^\rho\rightarrow {\mathcal H}_\nilp$
of \ref{section-preuve-2}. 
Let $q$ be the column of $T$ containing $n$.
Let us show that the image of $\Phi_T$ is the $Z(\nilp)$-orbit ${\mathcal H}_{\nilp,q}$. \\
We use the same notation as in \ref{section-preuve-2}.
A hyperplane $H$ in the image of $\Phi_T$
is such that the Young diagram $Y(\nilp_{|H})$ associated to the restriction $\nilp_{|H}$
is equal to $Y'$.
As $Y'$ is obtained from $Y(\nilp)$ by removing one box in the $q$-th column,
it follows
$\mathrm{ker}\,\nilp^{q-1}\subset H$ and $\mathrm{ker}\,\nilp^{q}\not\subset H$.
Thus $\Phi_T(\,\B_{\nilp,T}^\rho\,)\subset {\mathcal H}_{\nilp,q}$.
As ${\mathcal H}_{\nilp,q}$ is a $Z(\nilp)$-orbit and as
$\Phi_T$ is $Z(\nilp)$-equivariant, we get $\Phi_T(\,\B_{\nilp,T}^\rho\,)={\mathcal H}_{\nilp,q}$.

\paragraphe
\label{section-preuve-5}
Let $H'=H_{\mu_q}$ be one representative of the open cell ${\mathcal C}(\mu_q)\subset{\mathcal H}_{\nilp,q}$.
Let $\nilp'=\nilp_{|H'}$ be the restriction of $\nilp$.
Then $\Phi_T^{-1}(H')=\B_{\nilp',T'}^{\rho'}$.
By induction hypothesis $\B_{\nilp',T'}^{\rho'}$ is irreducible,
nonsingular and $\mathrm{dim}\,\B_{\nilp',T'}^{\rho'}=\mathrm{dim}\,\B_{\nilp'}$.
Moreover there is a cell decomposition $\B_{\nilp',T'}^{\rho'}=\bigsqcup_{\tau'}C^{\rho'}(\tau')$
parameterized by the row-standard tableaux $\tau'$ of shape $Y'$ 
with standardization $\mathrm{st}(\tau')=T'$, and such that $\mathrm{dim}\,C^{\rho'}(\tau')=\mathrm{dim}\,\B_{\nilp'}-n_\mathrm{inv}(\tau')$.

As $Z(\nilp)$ acts transitively on the image of $\Phi_T$,
it follows that the algebraic map 
\[\Xi:Z(\nilp)\times \B_{\nilp',T'}^{\rho'}\rightarrow \B_{\nilp,T}^\rho,\ 
(g,(V_0,...,V_{n-1})\,)\mapsto (gV_0,...,gV_{n-1},V)
\]
is surjective.
Moreover the restriction of $\Xi$ to $U(\mu_q)\times \B_{\nilp',T'}^{\rho'}$
is an isomorphism of algebraic varieties onto $\Phi_T^{-1}({\mathcal C}(\mu_q))$,
the inverse image of the open cell of ${\mathcal H}_{\nilp,q}$.
As $Z(\nilp)$ and $\B_{\nilp',T'}^{\rho'}$ are irreducible,
the surjectivity of $\Xi$ implies that $\B_{\nilp,T}^\rho$ is irreducible.
The subsets $g\Phi_T^{-1}({\mathcal C}(\mu_q))$ for $g\in Z(\nilp)$
form a covering of $\B_{\nilp,T}^\rho$ by nonsingular open subsets,
hence $\B_{\nilp,T}^\rho$ is nonsingular.
Moreover we have 
\[\mathrm{dim}\,\B_{\nilp,T}^\rho=\mathrm{dim}\,\Phi_T^{-1}({\mathcal C}(\mu_q))
=\mu_q-1+\mathrm{dim}\,\B_{\nilp'}=\mathrm{dim}\,\B_{\nilp}\]
(by the formula in \ref{section-theoreme-principal}). The proof of Theorem \ref{theoreme-Spaltenstein-general} is complete.

For $l=\mu_{q+1}+1,...,\mu_q$ we can find $g_l\in Z(\nilp)$ such that $H_l=g_lH'$.
Then the restriction of $\Xi$ to $U(l)g_l\times \B_{\nilp',T'}^{\rho'}$
is an isomorphism of algebraic varieties onto $\Phi_T^{-1}({\mathcal C}(l))$.
For $\tau$ row-standard with $\mathrm{st}(\tau)=T$,
the entry $n$ is in the $q$-th column of $\tau$, at the end of some row.
Thus there is $p\in\{\mu_{q+1}+1,...,\mu_q\}$ such that $n$ is at the end of the
$p$-th row of $n$. Let $\tau'$ be the row-standard tableau obtained from
$\tau$ by putting the $p$-th row at the place of the $\mu_q$-th row
and moving by one rank to the up each row among the $(p+1)$-th,...,$\mu_q$-th ones.
Then $n$ is at the same place in $\tau'$ and $T$, we denote by $\tau''$ the subtableau
of $\tau'$ obtained by deleting $n$.
This is a row-standard tableau of shape $Y'$ and standardization $\mathrm{st}(\tau')=T'$.
We define
\[C^\rho(\tau)=\Xi(\,U(l)g_l\times C^{\rho'}(\tau'')\,).\]
We get thus a partition $\B_{\nilp,T}^\rho=\bigsqcup_\tau C^\rho(\tau)$
parameterized by row-standard tableaux $\tau$ of standardization $\mathrm{st}(\tau)=T$.
This partition is a cell decomposition since it is the product of two cell decompositions.
It follows from the definition of inversions 
that $n_\mathrm{inv}(\tau)=n_\mathrm{inv}(\tau'')+(\mu_q-l)$.
We deduce
\[\mathrm{dim}\,C^\rho(\tau)\mathrm{dim}\,{\mathcal C}(l)+\mathrm{dim}\,C^{\rho'}(\tau'')
=l-1+\mathrm{dim}\,\B_{\nilp'}-n_\mathrm{inv}(\tau'')
=\mathrm{dim}\,\B_\nilp-n_\mathrm{inv}(\tau).\]
Therefore this cell decomposition satisfies the required properties.
The proof of Theorem \ref{theoreme-decomposition-cellulaire}
is complete.
\cqfd

\subsubsection*{Remark. Another cell decomposition of $\B_\nilp$}

\paragraphe
\label{section-intersection-Schubert}
The construction of our cell decomposition 
relies on the Schubert cell decomposition of the
Grassmannian of hyperplanes of ${\mathcal H}(V/\mathrm{Im}\,\nilp)$,
and an inductive argument.
A construction of a different cell decomposition relies on the Schubert cell
decomposition of the flag variety $\B$.
Recall that, if $B\subset GL(V)$ is a Borel subgroup,
then the $B$-orbits of $\B$
form a cell decomposition $\B=\bigsqcup_{\sigma\in S_n}S(\sigma)$
parameterized by the permutations, and the cells are called Schubert cells.
We show that the intersection of the Schubert cells with $\B_\nilp$ gives
a cell decomposition of $\B_\nilp$
provided that the Borel subgroup $B$ is well chosen.
Our proof is different than in \cite{Shimomura}.

\medskip

We consider a Jordan basis of $\nilp$. Recall that $\lambda_1\geq...\geq \lambda_r$
denote the lengths of the Jordan blocks of $\nilp$.
Let us index the basis $(e_{p,q})$ for $p=1,...,r$ and $q=1,...,\lambda_p$
so that $(e_{p,q})_{q=1,...,\lambda_p}$ is the subbasis corresponding to the
$p$-th Jordan block and we have
\[
\nilp(e_{p,1})=0\mbox{ \ and \ }\nilp(e_{p,q})=e_{p,q-1} \mbox{ for }q=2,...,\lambda_p.
\]
Such a pair $(p,q)$ with $1\leq q\leq \lambda_p$ forms the coordinates
of some box in the diagram $Y(\nilp)$. The Jordan basis is thus
indexed on the boxes of $Y(\nilp)$.

We associate a particular flag $\F_\tau\in\B_\nilp$ to each
row-standard tableau $\tau$ of shape $Y(\nilp)$. %The
%tableau $\tau$ is equivalent 
%to a maximal chain of nonordered partitions
%$0=\pi^{(0)}\subset \pi^{(1)}\subset ...\subset \pi^{(n-1)}\subset\pi^{(n)}=\lambda(\nilp)$
For $p=1,...,r$ and $i=1,...,n$ let
$\pi^{(i)}_p$ be the number of entries among $1,...,i$ in the $p$-th row of $\tau$.
%(see \ref{}).
For $i=1,...,n$ we define the subspace
\[V_i=\langle e_{p,q}:p=1,...,r,\ q=1,...,\pi^{(i)}_p\rangle.\]
It is immediate that this subspace is stable by $\nilp$.
Finally let $\F_\tau=(V_0,...,V_n)$. \\
The basis being considered, as above, as indexed on the boxes of the diagram $Y(\nilp)$,
we see that $V_i$ is generated by the vectors 
associated to the boxes of numbers $1,...,i$ in $\tau$. 

Let $H\subset GL(V)$ be the subgroup of diagonal automorphisms in the basis.
The flags $\F_\tau$ are exactly the elements of $\B_\nilp$ which are fixed by 
$H$ for its natural action on flags. However $H$ does not leave $\B_\nilp$ invariant.
We introduce a subtorus $H'\subset H$ with the same fixed points, which
leaves $\B_\nilp$ invariant. To do this, set $\epsilon_{p,q}=nq-p$.
For $t\in\C^*$
let $h_t\in GL(V)$ be defined by $h_t(e_{p,q})=t^{\epsilon_{p,q}}e_{p,q}$
for $p=1,...,r$ and $q=1,...,\lambda_p$.
Let $H'=(h_t)_{t\in\C^*}$ be the subtorus so-obtained.
The $\epsilon_{p,q}$'s are pairwise distinct (since $1\leq p\leq n$)
hence $H'$ has the same fixed points as $H$. 
Moreover we have $h_t^{-1}\nilp h_t=t^n\nilp$ for any $t$. 
As $t^n$ acts trivially on flags, it follows that $h_t$ leaves $\B_\nilp$
invariant. 

For any $\F\in\B_\nilp$,
as $\B_\nilp$ is a projective variety, the map $t\mapsto h_t\F$
admits a limit when $t\rightarrow \infty$, and this limit is a fixed point
for the action of $H'$.
For $\tau$ row-standard, write
\[S(\tau)=\{\F\in\B_\nilp:\mathrm{lim}_{t\rightarrow \infty}h_t\F=\F_\tau\}.\]
We get a partition $\B_\nilp=\bigsqcup_\tau S(\tau)$ parameterized by row-standard tableaux.

Write $\{(p,q):$ $p=1,...,r$, $q=1,...,\lambda_p\}=\{(p_i,q_i):i=1,...,n\}$
so that we have $\epsilon_{p_1,q_1}<...<\epsilon_{p_n,q_n}$
Write $e_i=e_{p_i,q_i}$. Let $B\subset GL(V)$ be the Borel subgroup
of lower triangular automorphisms in the basis $(e_1,...,e_n)$.
Then the set $S(\tau)$ in the partition is the intersection between
$\B_\nilp$ and the Schubert cell $B\F_\tau\subset \B$.

We see that the flag $\F_\tau$ belongs to the Spaltenstein subset $\B_\nilp^T$
for $T=\mathrm{st}(\tau)$.
Let $P=\{g\in GL(V):g(\mathrm{ker}\,\nilp^q)=\mathrm{ker}\,\nilp^q\}$. This is a parabolic
subgroup of $GL(V)$. We can see that each $\B_\nilp^T$ in the Spaltenstein
partition of $\B_\nilp$ is the intersection between $\B_\nilp$ and some $P$-orbit
of the flag variety. Observe that $B\subset P$. Then we obtain
$S(\tau)\subset \B_\nilp^T$.

In particular we have $S(\tau)=\{\F\in \B_\nilp^T:\mathrm{lim}_{t\rightarrow \infty}h_t\F=\F_\tau\}$. The subset $\B_\nilp^T\subset\B_\nilp$ is open and nonsingular.
By \cite[\S 4]{Birula} it follows that $S(\tau)$ is isomorphic to an affine space.

\medskip

Therefore the $S(\tau)$'s form a cell decomposition of $\B_\nilp$ parameterized
by row-standard tableaux. Moreover $\B_\nilp^T=\bigsqcup_{\tau}S(\tau)$ where the union
is taken over tableaux $\tau$ of standardization $\mathrm{st}(\tau)=T$.
This cell decomposition is different than the decomposition of Theorem
\ref{theoreme-decomposition-cellulaire}. Indeed for $\tau=\mbox{\scriptsize $\young(34,12)$}$
we see that
$S(\tau)=\{\F_\tau\}$ is a cell of codimension $2$,
whereas $n_{\mathrm{inv}}(\tau)=1$.
The dimension of cells in this decomposition
is given in \cite[\S 5.10]{Nakajima}.
%or \cite{LF} \S 3.

%\paragraphe
%In the previous paragraph we show that the flags $\F_\tau$'s
%which belong to $\B_\nilp^T$ correspond to the row-standard tableau
%$\tau$ such that $\mathrm{st}(\tau)=T$.
%More generally
%for $\rho\in{\mathcal R}_n$
%write $\rho=(i_k,j_k)$. For $k=1,...,n$ let $\rho'_k= i_{k-1}$ or $\rho'_k=j_k$
%depending on whether $i_k<i_{k-1}$ or $j_k>j_{k-1}$.
%Thus $\{\rho'_1,...,\rho'_k\}=\{i_k+1,...,j_k\}$ for any $k$.
%We define a row-standard tableau $\rho\star\tau$ in two steps.
%First we replace each entry $i$ by $\rho'_i$.
%Second we arrange each row in increasing order,
%and we denote by $\rho\star\tau$ the tableau so-obtained.
%Then we see that the flag $\F_{\rho\star\tau}$ belongs to the subset
%$\B_{\nilp,T}^\rho$ for $T=\mathrm{st}(\tau)$.

%The cell decomposition $\B_{\nilp,T}^\rho=\bigsqcup_{\tau}C^\rho(\tau)$
%of Theorem \ref{theoreme-decomposition-cellulaire}
%can be constructed so that
%in addition the cell $C^\rho(\tau)$ contains the flag $\F_{\rho\star\tau}$
%(see \cite{LF} Theorem 23).

\partie{Calculation of Betti numbers}

\paragraphe
\label{cohomology} 
Let $X$ be an algebraic variety.
We consider the classical cohomology of sheafs
(see \cite{Kashiwara} for example).
Let $H^*(X,\Q)$ denote the cohomology space with rational coefficients
and let $H_c^*(X,\Q)$ denote the rational cohomology with compact support
(both coincide when $X$ is projective).
The following proposition recalls that the knowledge of a cell
decomposition of $X$
allows to compute the Betti numbers of $X$
(see \cite[\S 4.6]{Kashiwara}).

\paragraph{Proposition}{\em
Let $X$ be an algebraic variety on $\C$ 
which admits a cell decomposition $X=\bigsqcup_{i\in I}Z_i$.
For $m\in\N$ let $r_m$ be the number of $m$-dimensional cells. \\
(a) We have $H_c^l(X,\Q)=0$ for $l$ odd and 
$\mathrm{dim}\, H_c^{2m}(X,\Q)=r_m$ for any $m\in\N$.\\
(b) If $X$ is projective, then we have $H^l(X,\Q)=0$ for $l$ odd
and $\mathrm{dim}\, H^{2m}(X,\Q)=r_m$ for any $m\in\N$.}

\paragraphe
\label{bettiBuTrho}
Let $d=\mathrm{dim}\,\B_\nilp$ (see \ref{section-theoreme-principal}).
We fix $\rho\in {\mathcal R}_n$ and $T$ a standard tableau.
By Theorem \ref{theoreme-decomposition-cellulaire} and Proposition \ref{cohomology},
for any $m\in\N$, we get the formula
\[\mathrm{dim}\, H_c^{2m}(\B_{\nilp,T}^\rho,\Q)=\#\{\tau\mbox{ row-standard}:\mathrm{st}(\tau)=T,\ n_\mathrm{inv}(\tau)=d-m\}.\]
Let $b_m^T=\mathrm{dim}\, H_c^{2m}(\B_{\nilp,T}^\rho,\Q)$.
For $i=1,...,n$, let $q_i$ be the number of the column of $T$ containing $i$
and let $p_i$ be the number of rows of length $q_i$
in the subtableau $T[1,...,i]$.
By Proposition \ref{distribution-longueur} we have
\[b_{d-m}^T=\#\left\{(\kappa_1,...,\kappa_n):0\leq\kappa_i\leq p_i-1,\ \kappa_1+...+\kappa_n=m\right\}
\quad
\forall m=0,...,d.\]
Let
$\chi^T(x)=\sum_{m=0}^d b_{d-m}^T\,x^m$.
For $p\in \N$ we write $[p]_x=1+x+...+x^{p-1}$.
We get:

\paragraph{Proposition}{\em
\label{proposition_2}
We have \
$\displaystyle{\chi^T(x)=\prod_{i=1}^n[p_i]_x}$.}

\paragraphe
We deduce the Betti numbers of certain components of $\B_\nilp$.
Let $T_\mathrm{min}$ be the 
minimal standard tableau of shape $Y(\nilp)$
for the dominance relation (see \ref{BuTrho-alphapartition}).
By Corollary \ref{theoreme-Spaltenstein-general}
the subset $\B_{\nilp,T_\mathrm{min}}^\rho\subset \B_\nilp$ is a nonsingular irreducible component.
The polynomial $\chi^{T_\mathrm{min}}(x)$ is its Poincar\'e polynomial.
Let $\mu_1,...,\mu_s$ be the lengths of the columns of $Y(\nilp)$.
Then $\chi^{T_\mathrm{min}}(x)$ is
\[\chi^{T_\mathrm{min}}(x)=\prod_{q=1}^s [\mu_q]_x!\]
where, for $m\in\N$, we write $[m]_x!=\prod_{p=1}^m[p]_x$.

\medskip
\noindent{\em Example.}
Suppose $Y(\nilp)=\mbox{\scriptsize $\yng(2,2,1)$}$, thus $T_\mathrm{min}=\mbox{\scriptsize $\young(14,25,3)$}$.
We get
$\chi^{T_\mathrm{min}}(x)={[2]_x}^2\cdot [3]_x=1+3x+4x^2+3x^3+x^4$.

\paragraphe
Let $d=\mathrm{dim}\,\B_\nilp$.
Let $b_m=\mathrm{dim}\, H^{2m}(\B_\nilp,\Q)$.
Set $\chi(x)=\sum_{m=0}^d b_{d-m}\,x^m$.
By Theorem \ref{theoreme-decomposition-cellulaire} and Proposition \ref{cohomology},
we get the following

\paragraph{Proposition} {\em
We have $\chi(x)=\sum_T\chi^T(x)$, where the sum is taken for $T$ running over the set of standard tableaux
of shape $Y(\nilp)$.}

\paragraphe {\em Inductive formula.}
If $Y$ is a Young diagram,
then we write $\chi(Y)(x):=\chi(x)=\sum_{m=0}^d\mathrm{dim}\, H^{2(d-m)}(\B_\nilp,\Q)x^m$ the above polynomial for $Y=Y(\nilp)$.
A box of $Y$ is said to be a corner it has no neighbor on the right or below.
Let $C(Y)\subset Y$ be the set of corners of $Y$.
Removing a corner $c$, we get a subdiagram $Y\setminus c\subset Y$.
For $c\in C(Y)$, let $q_c$ be the number of the column of $Y$ containing $c$
and let $p_c$ be the number of rows of $Y$ of length $q_c$.
We have the following inductive formula for the polynomial $\chi(Y)(x)$.

\paragraph{Proposition}{\em We have $\displaystyle{
\chi(Y)(x)=\sum_{c\in C(Y)}[p_c]_x\,\chi(Y\setminus c)(x)}$.}

\medskip
\noindent
{\em Proof.}
Let
$\chi(Y)_c(x)=\sum_{T\in \mathrm{Tab}_c(Y)}\chi^T(Y)(x)$,
where $\mathrm{Tab}_c(Y)$ denotes the set of standard tableaux of 
shape $Y$ such that $c$ contains the entry $n$.
We have thus $\chi(Y)(x)=\sum_{c\in C(Y)}\chi(Y)_c(x)$.
The set $\mathrm{Tab}_c(Y)$ is in bijection with the set of standard tableaux of shape $Y\setminus c$,
and, by Proposition \ref{proposition_2}, we see that $\chi(Y)_c(x)=[p_c]_x\,\chi(Y\setminus c)(x)$.
The proof is complete.
\cqfd

\bigskip
\noindent
{\em Example.}
We have computed $\chi(Y)(x)$ by induction, for several forms of $Y$:

\smallskip
\renewcommand{\arraystretch}{0.3}
{\scriptsize $\overbrace{\begin{array}{l}\young(\ \ )\cdots\young(\ )\end{array}}^n$}
$\displaystyle{\chi(Y)(x)=1}$
\qquad
{\scriptsize $\overbrace{\begin{array}{l}\young(\ \ )\cdots\young(\ )\\ \young(\ )\end{array}}^{s}$}
$\displaystyle{\chi(Y)(x)=s+x}$

\medskip

{\scriptsize $n\left\{\begin{array}{l}\young(\ ,\ )\\ { }^{\ \vdots}\\ \young(\ )\end{array}\right.$}
$\displaystyle{\chi(Y)(x)=[n]_x!}$
\quad
{\scriptsize $r\left\{\begin{array}{l}\young(\ \ ,\ )\\ { }^{\ \vdots}\\ \young(\ )\end{array}\right.$}
$\displaystyle{\chi(Y)(x)=[r-1]_x!\sum_{p=0}^{r-1}(r-p)x^p}$

\medskip
\noindent
and more generally:

\smallskip
{\scriptsize $r\left\{\begin{array}{l} \ \\ \ \\ \ \\ \ \\ \ \\ \ \\ \ \\ \ \\
\ \\ \ \\ \ \\ \ \\ \ \\ \ \\ \ \\ \ \\ \end{array}\right.\!\!\!\!\!\!\!\!\!\!\!
\overbrace{\begin{array}{l}\young(\ \ )\cdots\young(\ )\\ \young(\ ) \\ \ \vdots \\ \young(\ )\end{array}}^{s}$}
$\displaystyle{\chi(Y)(x)%=[r-1]_x!\,\sum_{p=1}^rC_{s+2+r-p}^{r-p}[p]_x
=[r-1]_x!\,\sum_{p=0}^{r-1}\left({ }^{s+p-2}_{\ \ \,p}\right)[r-p]_x}$
\quad (for $s\geq 2$)

\smallskip
{\scriptsize $\begin{array}{l}
\overbrace{\begin{array}{l}\young(\ \ )\cdots\young(\ \ )\end{array}}^{s} \\
\underbrace{\begin{array}{l}\young(\ )\cdots\young(\ \ )\end{array}}_{t}
\end{array}$}
$\chi(Y)(x)=[2]_x^t+{\displaystyle{\sum_{p=1}^t}}
\left({ }^{s+p-1}_{\ p-1}\right)\frac{s-p}{p}\,[2]_x^{t-p}$

\medskip
\noindent
and also:

\smallskip
{\scriptsize $\overbrace{\begin{array}{l}\young(\ \ \ )\cdots\young(\ )\\ \young(\ \ ,\ )\end{array}}^{s}$}
$\chi(Y)(x)\frac{s-3}{3}\left({ }^{s+2}_{\ \,2}\right)+
\frac{s+3}{3}\left({ }^{s-1}_{\ \,2}\right)[2]_x
+\left({ }^{s+1}_{\ \,2}\right)[2]_x^2+s[2]_x[3]_x$.

\end{document}